\documentclass[11pt]{article}

\usepackage[cp1250]{inputenc}
\usepackage{multirow}
\usepackage{amsthm}
\usepackage{amsmath,amsfonts,amssymb,graphics,verbatim,fullpage}

\newcommand{\cA}{\mathcal{A}}

\newcommand{\cF}{\mathcal{F}}
\newcommand{\cG}{\mathcal{G}}
\newcommand{\cH}{\mathcal{H}}

\newcommand{\IE}{\mathbb{E}}

\newcommand{\IF}{\mathbb{F}}
\newcommand{\IG}{\mathbb{G}}

\newcommand{\IH}{\mathbb{H}}

\newcommand{\IP}{\mathbb{P}}
\newcommand{\IQ}{\mathbb{Q}}

\newcommand{\R}{\mathbb{R}}

\newcommand{\be}{\begin{eqnarray*}}
\newcommand{\ee}{\end{eqnarray*}}
\newcommand{\ben}{\begin{eqnarray}}
\newcommand{\een}{\end{eqnarray}}

\newcommand{\halmos}{\vspace{3mm} \hfill $\Box$}

\newtheorem{theorem}{Theorem}[section]
\newtheorem{lemma}[theorem]{Lemma}

\theoremstyle{definition}

\newtheorem{problem}[theorem]{Problem}

\theoremstyle{remark}
\newtheorem{remark}[theorem]{Remark}

\numberwithin{equation}{section}

\begin{document}

\title{Quantile hedging for an insider}

\author{Przemys\l aw Klusik
\footnote{University of Wroc\l aw, pl.\ Grunwaldzki 2/4, 50-384 Wroc\l aw, Poland, E-mail: przemyslaw.klusik@math.uni.wroc.pl}
\qquad Zbigniew Palmowski\footnote{University of Wroc\l aw, pl.\ Grunwaldzki 2/4, 50-384 Wroc\l aw, Poland, E-mail: zbigniew.palmowski@math.uni.wroc.pl} \qquad Jakub Zwierz\footnote{University of Wroc\l aw, pl.\ Grunwaldzki 2/4, 50-384 Wroc\l aw, Poland, E-mail: jakub.zwierz@math.uni.wroc.pl}
}

\maketitle

\begin{abstract}
In this paper we consider the problem of the quantile hedging from the point of view
of a better informed agent acting on the market. The additional knowledge of the agent is modelled by a filtration initially enlarged by some random variable. By using equivalent martingale measures introduced in \cite{A00} and \cite{AIS98} we solve the problem for the complete case, by extending the results obtained in \cite{FL99} to the insider context. Finally, we consider the examples with the explicit calculations within the standard Black-Scholes model.

\medskip


\medskip

\end{abstract}

\section{Introduction}

A trader on the stock market is usually assumed to make
his decisions relying on all the information which is generated
by the market events.
However it is registered that some people have more
detailed information than others, in the sense that they
act with the present time knowledge of some future event.
This is the so-called insider information and those dealers
taking advantage of it are the insiders.
The financial markets with economic agents possessing additional knowledge have been studied in
a number of papers (see e.g. \cite{A00, AIS98, GP99, Z07}).
We take approach originated in \cite{DuffieHuang} and \cite{Pikovsky}
assuming that the insider possesses some extra information stored in the random variable $G$
known at the beginning of the trading interval and not available to the regular trader.
The typical examples of $G$ are $G=S_{T+\delta}$, $G=\mathbf{1}_{[a,b]}(S_{T+\delta})$ or
$G=\sup_{t\in[0,T+\delta]}S_t$ ($\delta >0$),
where $S$ is a semimartingale representing the discounted stock price process and $T$ is fixed time
horizon till which the insider is allowed to trade.

In this paper we show how much better and with which strategies an insider
can perform on the market if he uses optimally the extra
information he has at his disposal.
The problem of pricing and perfect hedging of contingent claims is well understood in the
context of arbitrage-free models which are complete. In such models every contingent
claim can be replicated by a self-financing trading strategy. The cost
of replication equals the discounted expectation of
the claim under the unique equivalent martingale measure.
Moreover, this cost is the same for the insider and the regular trader.
Therefore instead of this strategy we will employ the quantile
hedging strategy of an insider for the replication, following idea of
F\"{o}llmer and Leukert \cite{FL99, FL00}.
That is, we will seek for the self-financing strategy that
\begin{itemize}
\item maximizes the probability of success of hedge under a given initial capital or
\item minimizes the initial capital under a given lower bound of the probability of the successful hedge.
\end{itemize}

This is the case when the insider is unwilling to put up the initial amount of capital required by a
perfect hedging.
This approach might be also seen as a
dynamic version of the VaR.

We use powerful technique of {\it grossissement de filtrations}
developed by Yor, Jeulin and Jacod \cite{J80, JY85}
and utilize the results of Amendinger
\cite{A00} and Amendinger et al. \cite{AIS98}.

The paper is organized as follows. In Section \ref{market} we present the main results. In Section \ref{numerical}
we analyze in detail some examples. Finally, in Section \ref{proofs} we give the proofs of the main results.

\section{Main results}\label{market}

Let $(\Omega,\cF,\IF,\IP)$ be a complete probability space
and $S=(S_t)_{t\geq 0}$ be an $(\IF,\IP)$-semimartingale representing the discounted stock price process.
Assume that the filtration $\IF=(\cF_t)_{t\geq 0}$ is the natural filtration of $S$ satisfying usual conditions with the trivial $\sigma$-algebra $\cF_0$. Thus, the regular trader makes his portfolio decisions according to the information flow
$\IF$. In addition to the regular trader we will consider the insider, whose knowledge will be modelled by the \textit{initial enlargement} of $\IF$, that is filtration $\IG=(\cG_t)_{t\geq 0}$ given by:
$$\cG_t=\cF_t \vee
\sigma (G),
$$
where $G$ is an $\cF$-measurable random variable. In particular, $G$ can be an
$\cF_{T+\delta}$-measurable random variable ($\delta >0$) for $T$ being a fixed time horizon representing
the expiry date of the hedged contingent claim.

We will assume that the market is complete and arbitrage-free for the regular trader, hence there exists a unique equivalent martingale measure $\IQ_{\IF}$ such that $S$ is an $(\IF,\IQ_{\IF})$-martingale on $[0,T]$.
Denote by $(Z_{t}^{\IF})_{t\in[0,T]}$ the density process of $\IQ_\IF$ with respect to $\IP$, i.e.:
$$Z_{t}^{\IF}=\left.\frac{{\rm d}\IQ_\IF}{{\rm d}\IP}\right|_{\cF_t}.$$

We will consider the contingent claim $H$ being $\cF_T$-measurable, nonnegative random variable and the replicating
investment strategies for insider, which are expressed in terms of the integrals with respect to $S$.
To define them properly we assume that $S$ is a $(\IG,\IP)$-semimartingale
which follows from the requirement:
\begin{equation}
\IP(G\in\cdot|\cF_t)\ll \IP(G\in\cdot)\qquad \IP-{\rm a.s.} \label{abs_cont}
\end{equation}
for all $t\in [0,T]$ (see e.g. \cite{A00} and \cite{JY85}).
In fact, we assume from now on more, that is that the measure $\IP(G\in\cdot|\cF_t)$ and the law of $G$ are equivalent
for all $t\in [0,T]$:
\begin{equation}
\IP(G\in\cdot|\cF_t)\sim \IP(G\in\cdot)\qquad \IP-{\rm a.s.} \label{equiv}
\end{equation}
Under the condition (\ref{equiv}) there exists equivalent
$\IG$-martinagle measure $\IQ_{\IG}$ defined by:
\begin{equation}\label{QG}
\IQ_\IG(A)=\int_{A}\frac{Z_{T}^{\IF}}{p_{T}^{G}}(\omega){\rm d}\IP(\omega),\qquad A\in\cG_T,\end{equation}
where
$$p_{t}^{x}\;\IP(G\in {\rm d}x)$$
is a version of $\IP(G\in {\rm d}x|\cF_t);$
see \cite{A00} and \cite{AIS98} (and also Theorems \ref{ther1} and \ref{ther2}).

For $\IH\in\{\IF,\IG\}$ we will consider only self-financing admissible trading strategies $(V_0,\xi)$ on $[0,T]$
for which the value process
$$V_{t}=V_{0}+\int_{0}^{t}\xi_u\;{\rm d}S_u,\qquad t\in[0,T],$$
is well defined, where an initial capital $V_{0}\geq 0$ is $\cH_0$-measurable, a process $\xi$ is $\IH$-predictable
and
$$V_t\geq 0,\qquad \IP-{\rm a.s.}$$
for all $t\in [0,T]$.
Denote all admissible strategies associated to the filtration $\IH\in\{\IF,\IG\}$ by $\cA^{\IH}$.

Under assumption (\ref{equiv}) the insider can perfectly replicate the contingent
claim $H\in L^{1}(\IQ_\IF)\cap L^{1}(\IQ_\IG)$:
$$\IE_{\IQ_\IG}(H|\cG_t)=H_{0}+\int_{0}^{t}\xi_{u}\;{\rm d}S_{u},\qquad\IP-{\rm a.s.}$$
where $H_{0}=\IE_{\IQ_\IG}(H|\cG_0)$. Moreover, from \cite{A00} and \cite{AIS98} it follows that $H_{0}=\IE_{\IQ_\IG}H$ (see Theorem \ref{ther1}).
In this paper we will analyze the case when the insider is unwilling to pay the initial capital
$H_0$ required by a perfect hedge.
We will consider the following pair of dual problems.

\begin{problem}\label{P1}
Let $\alpha$ be a given $\cG_0$-measurable random variable taking values in $[0,1]$. We are looking for a strategy $(\alpha\IE_{\IQ_\IG}H,\xi)\in\cA^\IG$ which maximizes for any realization of $G$ the insider's probability of a successful hedge
$$\IP\left(\left.\alpha\IE_{\IQ_\IG}H+\int_0^T\xi_t \;{\rm d}S_t\geq H\right|\cG_0\right),\qquad \IP-{\rm a.s.}$$
\end{problem}
\begin{problem}\label{P2}
Let $\epsilon$ be a given $\cG_0$-measurable random variable taking values in $[0,1]$. We are looking for a minimal $\cG_0$-measurable random variable $\alpha$ for which there exists $\xi$ such that $(\alpha\IE_{\IQ_\IG}H,\xi)\in\cA^\IG$ and
\begin{equation}\label{P2inequality}\IP\left(\left.\alpha\IE_{\IQ_\IG}H+\int_0^T\xi_t \;{\rm d}S_t\geq H \right|\cG_0\right)\geq 1-\epsilon,\qquad \IP-{\rm a.s.}\end{equation}
\end{problem}
\begin{remark}
Recall that in the quantile hedging problem for the usual trader we maximize the objective probability $\IP(\alpha\IE_{\IQ_\IF}H+\int_0^T\xi_t \;{\rm d}S_t\geq H)$, where $\alpha$ is number from $[0,1]$.
In the Problems \ref{P1}-\ref{P2} we use conditional probability, since
now the insider's perception of the market at time $t=0$ depends on the knowledge described by $\cG_0$.
\end{remark}
The set
\begin{equation}\label{succset}\left\{\alpha\IE_{\IQ_\IG}H+\int_0^T\xi_t \;{\rm d}S_t\geq H\right\}\end{equation}
we will call the \textit{success set}.
Denote
$$\IQ^*(A)=\frac{\IE_{\IQ_\IG}(H\mathbf{1}_A)}{\IE_{\IQ_\IG}(H)},\qquad A\in\cG_T.$$

The following theorems solve Problems \ref{P1} and \ref{P2}.
\begin{theorem}\label{maxprobab}
Suppose that there exists a $\cG_0$-measurable random variable $k$ such that
\begin{equation}\label{eqk}\IQ^*\left(\left.\left.\frac{\;{\rm d}\IQ^*}{\;{\rm d}\IP}\right|_{\cG_T}\leq k\right|\cG_0\right)=\alpha.\end{equation}
Then the maximal probability of a success set solving Problem \ref{P1} equals:
$$\IP\left(\left.\left.\frac{\;{\rm d}\IQ^*}{\;{\rm d}\IP}\right|_{\cG_T}\leq k\right|\cG_0\right)$$
and it is realized by the strategy
$$\left(\IE_{\IQ_\IG}\left[\left.H\mathbf{1}_{\left\{\left.\frac{\;{\rm d}\IQ^*}{\;{\rm d}\IP}\right|_{\cG_T}\leq k\right\}}\right|\cG_0\right],\tilde{\xi}\right),$$
which replicates the payoff $H\mathbf{1}_{\left\{\left.\frac{\;{\rm d}\IQ^*}{\;{\rm d}\IP}\right|_{\cG_T}\leq k\right\}}$.
\end{theorem}

\begin{theorem}\label{minCostQH}
Suppose that there exists a $\cG_0$-measurable random variable $k$ such that
\begin{equation}\label{eqk2}\IP\left(\left.\left.\frac{\;{\rm d}\IQ^*}{\;{\rm d}\IP}\right|_{\cG_T}\leq k\right|\cG_0\right)=1-\epsilon.\end{equation}
Then the minimal $\cG_0$-measurable random variable $\alpha$ solving Problem \ref{P2} equals:
$$\IQ^*\left(\left.\left.\frac{\;{\rm d}\IQ^*}{\;{\rm d}\IP}\right|_{\cG_T}\leq k\right|\cG_0\right)$$
and it is realized by the strategy
$$\left(
\IE_{\IQ_\IG}\left[\left.H\mathbf{1}_{\left\{\left.\frac{\;{\rm d}\IQ^*}{\;{\rm d}\IP}\right|_{\cG_T}\leq k\right\}}\right|\cG_0\right],\tilde{\xi}\right)$$ being the perfect hedge of $H\mathbf{1}_{\left\{\left.\frac{\;{\rm d}\IQ^*}{\;{\rm d}\IP}\right|_{\cG_T}\leq k\right\}}.$
\end{theorem}

\begin{remark}
The assumptions that there exists $k$ satisfying (\ref{eqk}) and (\ref{eqk2}) are satisfied if
$\IP(Z_T^{\IF}H=0|\cG_0)<\alpha$ and $\IP(Z_T^{\IF}H=0|\cG_0)<1-\epsilon$ respectively, and
$Z_T^{\IF}H$ has the conditional density on $\R_+$ given $G=g$ - see Section \ref{numerical} for the examples.
\end{remark}

The proofs of these theorems are given in Section \ref{proofs}.

\section{Numerical examples}\label{numerical}

In this section we consider the standard Black-Scholes model
in which the price evolution is described by the equation
$${\rm d}S_{t}=\sigma S_{t}{\rm d}W_{t}+\mu S_{t}{\rm d}t,$$
where $W$ is a Brownian motion, $\sigma,\mu>0$.
For simplicity we assume that interest rate is zero.
We analyze the Problem \ref{P2} for two examples of the insider information
and provide numerical results for pricing the vanilla call option, where
$$H=(S_T-K)^+$$
and $K$ is a strike price.

\subsection{The case of $G=W_{T+\delta}$}
It means that insider knows the stock price $G=S_{T+\delta}$
after the expiry date $T$.
In this case we have:
\begin{eqnarray*}
\lefteqn{\IP(W_{T+\delta}\in {\rm d}z|\cF_t)=\IP(W_{T+\delta}-W_{t}+W_{t}\in {\rm d}z|\cF_t)}\\
&&=  \frac{1}{\sqrt{2\pi(T+\delta-t)}}\exp\Big(-\frac{(z-W_{t})^{2}}{2(T+\delta-t)}\Big){\rm d}z\\
&& =  p_{t}^{z}\IP(W_{T+\delta}\in {\rm d}z),
\end{eqnarray*}
where
$$p_{t}^{z}=\sqrt{\frac{T+\delta}{T+\delta-t}}
\exp\Big(-\frac{(z-W_{t})^{2}}{2(T+\delta-t)}+\frac{z^{2}}{2(T+\delta)}\Big).$$
Therefore:
\begin{eqnarray*}
\lefteqn{ \left.\frac{{\rm d}\IQ_\IG}{{\rm d}\IP}\right|_{\cG_T}=\frac{Z_T^{\IF}}{p_T^G}=
\frac{\left.\frac{{\rm d}\IQ_\IF}{{\rm d}\IP}\right|_{\cF_T}}{p_T^G}=
\frac{\exp\Big(-\frac{\mu}{\sigma}W_T-\frac{1}{2}\big(\frac{\mu}{\sigma}\big)^{2}T\Big)
}{\sqrt{\frac{T+\delta}{\delta}}\exp\Big(-\frac{(W_{T+\delta}-W_T)^{2}}{2\delta}
+\frac{W_{T+\delta}^{2}}{2(T+\delta)}\Big)}}\\
&& =  \sqrt{\frac{\delta}{T+\delta}}\exp\Big(-\frac{\mu}{\sigma}W_T-\frac{1}{2}
\Big(\frac{\mu}{\sigma}\Big)^{2}T+\frac{(W_{T+\delta}-W_T)^{2}}{2\delta}-\frac{W_{T+\delta}^{2}}{2(T+\delta)}\Big)
\end{eqnarray*}
and
\begin{equation}\label{eq1ex1}\left.\frac{{\rm d}\IQ^*}{{\rm d}\IP}\right|_{\cG_T}
=\frac{H}{\IE_{\IQ_\IG}H}\left.\frac{{\rm d}\IQ_\IG}{{\rm d}\IP}\right|_{\cG_T}.\end{equation}
Note that 
$\left.\frac{{\rm d}\IQ^*}{{\rm d}\IP}\right|_{\cG_T}$ is a random variable with the conditional density
on $\R_+$ given $G=g$ and for given $\epsilon\in[0,1]$ we can find a $\cG_0$-measurable random variable $k$ such that
$\IP\left(\left.\left.\frac{{\rm d}\IQ^*}{{\rm d}\IP}\right|_{\cG_T}\leq k\right|\cG_0\right)=1-\epsilon$
if $\IP(H=0|\cG_0)<1-\epsilon$. Therefore, by Theorem \ref{minCostQH} the cost of the quantile hedging for the insider can be reduced in this case by the factor:
\begin{equation}\label{eq2ex1}\alpha=\IQ^*\left(\left.\left.\frac{\;{\rm d}\IQ^*}{\;{\rm d}\IP}\right|_{\cG_T}\leq k\right|\cG_0\right).\end{equation}
Below, we provide the values of $\alpha$ for $\mu=0.08$, $\sigma=0.25$, $S_0=100$, $K=110$, $T=0.25$, $\delta=0.02$, and different values of $G$ and $\epsilon$. In the programme we use simple fact that $\IE[f(W_T)|G=g]=f(g-W(\delta))$ for a measurable function $f$.\newline
\begin{center}
\begin{tabular}{|c|c|c|c|c|c|c|c|c|c|c|c|c|}
\hline
~ & \multicolumn{12}{|c|}{$G$}\\
\hline
~ & ~ & 105 & 106 & 107 & 108 & 109 & 110 & 111 & 112 & 113 & 114 & 115\\
\hline
\multirow{6}{*}{$\epsilon$} & 0.01 & 0.05 & 0.09 & 0.13 & 0.17 & 0.22 & 0.27 & 0.32 & 0.37 & 0.42 & 0.46 & 0.51\\
& 0.05 & $<0.01$ & 0.01 & 0.04 & 0.07 & 0.10 & 0.14 & 0.18 & 0.23 & 0.28 & 0.33 & 0.38\\
& 0.10 & $<0.01$ & $<0.01$ & 0.01 & 0.03 & 0.05 & 0.08 & 0.12 & 0.16 & 0.21 & 0.25 & 0.30\\
& 0.15 & $<0.01$ & $<0.01$ & $<0.01$ & 0.01 & 0.03 & 0.05 & 0.08 & 0.12 & 0.16 & 0.21 & 0.25\\
& 0.20 & $<0.01$ & $<0.01$ & $<0.01$ & $<0.01$ & 0.01 & 0.03 & 0.06 & 0.09 & 0.13 & 0.17 & 0.21\\
& 0.25 & $<0.01$ & $<0.01$ & $<0.01$ & $<0.01$ & $<0.01$ & 0.02 & 0.04 & 0.07 & 0.10 & 0.14 & 0.18\\
\hline
\end{tabular}
\end{center}
~\newline
\subsection{The case of $G=\mathbf{1}_{\left\{W_{T+\delta}\in[a,b]\right\}}$}

In this example the insider knows the range of the stock price $S_{T+\delta}$
after the expiry date $T$.
The straightforward calculation yields:
\begin{eqnarray*}
\IP(G=1|\cF_t) & = & \frac{1}{\sqrt{2\pi(T+\delta-t)}}\int_{a-W_t}^{b-W_t}\exp\Big(-\frac{u^2}{2(T+\delta-t)}\Big)\;{\rm d}u\\
~ & = & \Phi\left(\frac{b-W_t}{\sqrt{T+\delta-t}}\right)-\Phi\left(\frac{a-W_t}{\sqrt{T+\delta-t}}\right),\\
\end{eqnarray*}
where $\Phi$ is c.d.f. of the standard normal distribution. Thus,
$$p_t^1=\frac{\IP(G=1|\cF_t)}{\IP(G=1)}=\frac{\Phi\left(\frac{b-W_t}{\sqrt{T+\delta-t}}\right)-\Phi\left(\frac{a-W_t}{\sqrt{T+\delta-t}}\right)}{\Phi\left(\frac{b}{\sqrt{T+\delta}}\right)-\Phi\left(\frac{a}{\sqrt{T+\delta}}\right)},$$
and similarly
$$p_t^0=\frac{1+\Phi\left(\frac{a-W_t}{\sqrt{T+\delta-t}}\right)-\Phi\left(\frac{b-W_t}{\sqrt{T+\delta-t}}\right)}{1+\Phi\left(\frac{a}{\sqrt{T+\delta}}\right)-\Phi\left(\frac{b}{\sqrt{T+\delta}}\right)}.$$
Hence
\begin{eqnarray*}
\lefteqn{\left.\frac{{\rm d}\IQ_\IG}{{\rm d}\IP}\right|_{\cG_T} =  \frac{\exp\left(-\frac{\mu}{\sigma}W_T-\frac{1}{2}\left(\frac{\mu}{\sigma}\right)^2T\right)}{
p_T^0\mathbf{1}_{\{G=0\}}+p_T^1\mathbf{1}_{\{G=1\}}}}\\
&& =  \left.\exp\left(-\frac{\mu}{\sigma}W_T-\frac{1}{2}\left(\frac{\mu}{\sigma}\right)^2T\right)\right/
\Bigg(\mathbf{1}_{\{G=0\}}\frac{1+\Phi\left(\frac{a-W_T}{\sqrt{T+\delta-t}}\right)-\Phi\left(\frac{b-W_T}{
\sqrt{T+\delta-t}}\right)}{1+\Phi\left(\frac{a}{\sqrt{T+\delta}}\right)-\Phi\left(\frac{b}{\sqrt{T+\delta}}\right)}\\
&&\hspace{7cm} +\mathbf{1}_{\{G=1\}}\frac{\Phi\left(\frac{b-W_T}{\sqrt{T+\delta-t}}\right)-\Phi\left(\frac{a-W_T}{\sqrt{T+\delta-t}}\right)}{
\Phi\left(\frac{b}{\sqrt{T+\delta}}\right)-\Phi\left(\frac{a}{\sqrt{T+\delta}}\right)}\Bigg)
\end{eqnarray*}
and $\IQ^*$ and $\alpha$ are defined in (\ref{eq1ex1})-(\ref{eq2ex1}).
The table below provides the values of the optimal $\alpha$ for $\mu=0.08$, $\sigma=0.25$, $S_0=100$, $K=110$, $T=0.25$, $\delta=0.02$, $G=1$ and different values of $\epsilon$ and endpoints of interval $[a,b]$ for $S_{T+\delta}$. In the programme we use simple fact that $\IE[f(W_T)|G=1]=\IE[f(W_T), W_{T+\delta}\in[a,b]]/\IP(G=1)$ for a measurable function $f$
and to simulate numerator we choose only those trajectories for which $W_{T+\delta}\in[a,b]$.
\newline
\begin{center}
\begin{tabular}{|c|c|c|c|c|c|c|}
\hline
~ & \multicolumn{6}{|c|}{$[a,b]$}\\
\hline
~ & ~ & [109,111] & [108,112] & [107,113] & [112,114] & [106,108] \\
\hline
\multirow{6}{*}{$\epsilon$} & 0.01 & 0.272 & 0.284 & 0.296 & 0.413 & 0.135 \\
& 0.05 & 0.142 & 0.150 & 0.157 & 0.277 & 0.039 \\
& 0.10 & 0.087 & 0.088 & 0.095 & 0.209 & 0.010 \\
& 0.15 & 0.053 & 0.055 & 0.059 & 0.164 & 0.001 \\
& 0.20 & 0.032 & 0.033 & 0.034 & 0.129 & $<0.001$ \\
& 0.25 & 0.017 & 0.019 & 0.020 & 0.102 & $<0.001$ \\
\hline
\end{tabular}
\end{center}

\section{Proofs}\label{proofs}

Before we give the proofs of the main Theorems \ref{maxprobab} and \ref{minCostQH} we present few
introductory lemmas and theorems. We start with the result of \cite{A00} and \cite{AIS98} concerning the properties of the equivalent martingale measure $\IQ_\IG$ for the insider. We recall that we assume that
the condition (\ref{equiv}) is satisfied.
\begin{theorem}\label{ther1}
\begin{itemize}
\item[(i)] Process $Z_t^{\IG}:=\frac{Z^{\IF}_t}{p_t^G}$ is a $(\IG,\IP)$-martingale.
\item[(ii)] The measure $\IQ_\IG$ defined in (\ref{QG}) has the following properties:
\subitem(a) $\cF_T$ and $\sigma(G)$ are independent under $\IQ_\IG$;
\subitem(b) $\IQ_\IG=\IQ_\IF$ on $(\Omega,\cF_T)$ and $\IQ_\IG=\IP$ on $(\Omega,\sigma(G))$.
\end{itemize}
\end{theorem}
We are now in a position to state the theorem
which relates the martingale measures of the insider and the regular trader.
\begin{theorem}\label{ther2}
Let $X=(X_{t})_{t\geq 0}$ be an $\IF$-adapted process. The following statements are equivalent:
\begin{itemize}
\item[(i)]$X$ is an $(\IF,\IQ_\IF)$-martingale;
\item[(ii)]$X$ is an $(\IF,\IQ_\IG)$-martingale;
\item[(iii)]$X$ is a $(\IG,\IQ_\IG)$-martingale.
\end{itemize}
\end{theorem}
\begin{proof}
Equivalence (i) and (ii) follows from the fact that
$\IQ_\IF=\IQ_\IG$ on $\cF_T$. The implication $(iii)\Rightarrow (ii)$ is a consequence of the
tower property of the conditional expectation.
Finally, taking $A=A_{s}\cap\{\omega\in \Omega: G(\omega)\in B\}$ ($A_{s}\in\cF_s$, $B$ - Borel set), the
implication $(ii)\Rightarrow (iii)$ follows from the standard monotone class arguments and following equalities:
\begin{eqnarray*}
\lefteqn{\IE_{\IQ_\IG}(\mathbf{1}_{A}X_{t})=\IE_{\IQ_\IG}(\mathbf{1}_{A_s}\mathbf{1}_{\{G\in B\}}X_{t})=
\IQ_\IG(G\in B)
\IE_{\IQ_\IG}(\mathbf{1}_{A_s}X_{t})}
\\&&=\IQ_\IG(G\in B)
\IE_{\IQ_\IG}(\mathbf{1}_{A_s}X_{s})=
\IE_{\IQ_\IG}(\mathbf{1}_{A}X_{s}),\qquad s\leq t,
\end{eqnarray*}
where in the second equality we use Theorem \ref{ther1}(ii).
\end{proof}
\begin{remark}
We need the Theorem \ref{ther2} to guarantee the martingale representation for the insider's replicating strategy. Moreover, from this representation and the Theorem \ref{ther1} we can deduce that the cost of perfect
hedging for the insider is the same as for the regular trader, that is $\IE_{\IQ_\IG}[H|\cG_0]=\IE_{\IQ_\IF}H$.
\end{remark}
\begin{remark}
In general the Theorem \ref{ther2} is not true for the local martingales, since a localizing sequence $(\tau_{n})$ of $\IG$-stopping times is not a sequence of $\IF$-stopping times.
\end{remark}

\begin{lemma}\label{NPL}
Let $k$ be a positive $\cG_t$-measurable random variable. For every $A\in\cF$ such that
$$\IQ^*\left(\left.A\right|\cG_t\right)\leq\IQ^*
\left(\left.\left.\frac{{\rm d}\IQ^*}{{\rm d}\IP}\right|_{\cG_T}\leq k\right|\cG_t\right),\qquad \IP-{\rm a.s.}$$
we have:
$$\IP\left(\left.A\right|\cG_t\right)\leq\IP
\left(\left.\left.\frac{{\rm d}\IQ^*}{{\rm d}\IP}\right|_{\cG_T}\leq k\right|\cG_t\right),\qquad \IP-{\rm a.s.}$$
Similarly, if
$$\IP\left(\left.A\right|\cG_t\right)\geq\IP
\left(\left.\left.\frac{{\rm d}\IQ^*}{{\rm d}\IP}\right|_{\cG_T}\leq k\right|\cG_t\right),\qquad \IP-{\rm a.s.}$$
then
$$\IQ^*\left(\left.A\right|\cG_t\right)\geq\IQ^*
\left(\left.\left.\frac{{\rm d}\IQ^*}{{\rm d}\IP}\right|_{\cG_T}\leq k\right|\cG_t\right),\qquad \IP-{\rm a.s.}$$
\end{lemma}
\begin{proof}
Denote $\tilde{A}:=\left\{\left.\frac{{\rm d}\IQ^*}{{\rm d}\IP}\right|_{\cG_T}\leq k
\right\}$. Note that $$(\mathbf{1}_{\tilde{A}}-\mathbf{1}_A)
\left(k-\left.\frac{{\rm d}\IQ^*}{{\rm d}\IP}\right|_{\cG_T}\right)
\geq 0.$$ Thus
\begin{eqnarray*}
\left.\frac{\;{\rm d}\IQ^*}{\;{\rm d}\IP}\right|_{\cG_t}\left(\IQ^*(\tilde{A}|\cG_t)-\IQ^*(A|\cG_t)\right)
&=&\left.\frac{\;{\rm d}\IQ^*}{\;{\rm d}\IP}\right|_{\cG_t}\IE_{\IQ^*}\left((\mathbf{1}_{\tilde{A}}-\mathbf{1}_{A})\left.\right|\cG_t\right)\\
&=&\IE_\IP\left(\left.\left.\frac{{\rm d}\IQ^*}{{\rm d}\IP}\right|_{\cG_T}(\mathbf{1}_{\tilde{A}}-\mathbf{1}_{A})\right|\cG_t\right)\\
&\leq&k\left(\IP(\tilde{A}|\cG_t)-\IP(A|\cG_t)\right),
\end{eqnarray*}
which completes the proof.
\end{proof}

\begin{lemma}\label{equalone}
The following holds true:
$$\left.\frac{{\rm d}\IQ^*}{{\rm d}\IQ_\IG}\right|_{\cG_0}=1.$$
\end{lemma}
\begin{proof}
Note that for $A=\{\omega\in\Omega: G\in B\}\in \cG_0$ ($B$ is a Borel set) we have
\begin{eqnarray*}
\lefteqn{\IQ^*(A)=\IE_{\IQ_{\IG}}\left[\frac{{\rm d}\IQ^*}{{\rm d}\IQ_\IG}\mathbf{1}_{A}\right]=\frac{\IE_{\IQ_{\IG}}\left[H\mathbf{1}_{A}\right]}{\IE_{\IQ_{\IG}}H}}\\
&&=
\frac{\IE_{\IQ_{\IG}}\left[H\right]\IQ_{\IG}(A)}{\IE_{\IQ_{\IG}}H}=\IQ_{\IG}(A),
\end{eqnarray*}
where in the last but one equality we use Theorem \ref{ther1}.
\end{proof}

{\it Proof of Theorem \ref{maxprobab}}

Consider the value process $V_{t}=\alpha\IE_{\IQ_\IG}H+\int_{0}^{t}\xi_{u}\;{\rm d}S_{u}$ for
any strategy $(\alpha\IE_{\IQ_\IG}H,\xi)\in\cA^\IG$. Note that for its success set $A$ defined in (\ref{succset})
we have:
$$V_T\geq H\mathbf{1}_A.$$
Moreover, by the Theorem \ref{ther2} the process $V_{t}$ is a nonnegative $(\IG,\IQ_\IG)$-local martingale, hence it is a $(\IG,\IQ_\IG)$-supermartingale and
$$\alpha\IE_{\IQ_\IG}H=\alpha V_0\geq \IE_{\IQ_\IG}(V_T|\cG_0)\geq \IE_{\IQ_\IG}(H\mathbf{1}_{A}|\cG_0).$$
Thus, from Lemmas \ref{equalone} and \ref{NPL},
$$\IQ^*(A|\cG_0)\leq\frac{\alpha}{\left.\frac{{\rm d}\IQ^*}{{\rm d}\IQ_\IG}\right|_{\cG_0}}=\alpha,\qquad \IP-{\rm a.s.}$$
and therefore
\ben
\IP(A|\cG_0)\leq\IP\left(\left.\left.\frac{{\rm d}\IQ^*}{{\rm d}\IP}\right|_{\cG_T}\leq k\right|\cG_0\right),\qquad \IP-{\rm a.s.} \label{P_ineq}
\een
It remains to show that $\left(\IE_{\IQ_\IG}\left[H \left.\mathbf{1}_{\left\{\left.\frac{{\rm d}\IQ^*}{{\rm d}\IP}\right|_{\cG_T}\leq k\right\}}\right|\cG_0\right],\tilde{\xi}\right)\in\cA^\IG$, and that this strategy attains the upper bound
(\ref{P_ineq}). The first statement
follows directly from the definition of $\tilde{\xi}$:
$$
\IE_{\IQ_\IG}\left(H \left.\mathbf{1}_{\left\{\left.\frac{{\rm d}\IQ^*}{{\rm d}\IP}\right|_{\cG_T}\leq k\right\}}\right|\cG_0\right)
+\int_0^t\tilde{\xi}_u\;{\rm d}S_u=\IE_{\IQ_\IG}\left(H \left.\mathbf{1}_{\left\{\left.\frac{{\rm d}\IQ^*}{{\rm d}\IP}\right|_{\cG_T}\leq k\right\}}\right|\cG_t\right)\geq 0.$$
Moreover,
\ben
\lefteqn{\IP\left(\left.\IE_{\IQ_\IG}\left[H \left.\mathbf{1}_{\left\{\left.\frac{{\rm d}\IQ^*}{{\rm d}\IP}\right|_{\cG_T}\leq k\right\}}\right|\cG_0\right]
+\int_0^T\tilde{\xi}_u\;{\rm d}S_u\geq H\right|\cG_0\right)}\nonumber\\&&=
\IP\left(\left.H \mathbf{1}_{\left\{\left.\frac{{\rm d}\IQ^*}{{\rm d}\IP}\right|_{\cG_T}\leq k\right\}}\geq H\right|\cG_0\right)\geq
\IP\left(\left. \left.\frac{{\rm d}\IQ^*}{{\rm d}\IP}\right|_{\cG_T}\leq k\right|\cG_0\right),\nonumber
\een
which completes the proof in view of (\ref{P_ineq}).
\halmos

{\it Proof of Theorem \ref{minCostQH}}

Observe that for any $(\alpha\IE_{\IQ_\IG}H,\xi)\in\cA^\IG$ we
have: \ben
\lefteqn{\IQ^*\left(\left.\alpha\IE_{\IQ_\IG}H+\int_0^T\xi_u\;{\rm
d}S_u\geq H\right|\cG_0\right)}\nonumber\\&&
=\frac{1}{\left.\frac{{\rm d}\IQ^*}{{\rm
d}\IQ_\IG}\right|_{\cG_0}}
\IE_{\IQ_\IG}\left(\left.\left.\frac{{\rm d}\IQ^*}{{\rm
d}\IQ_\IG}\right|_{\cG_T}\mathbf{1}_{\left\{\alpha\IE_{\IQ_\IG}H+\int_0^T\xi_u\;{\rm
d}S_u\geq H\right\}}\right|\cG_0\right)\nonumber\\&&
\leq\frac{\IE_{\IQ_\IG}\left(\left.\alpha\IE_{\IQ_\IG}H+\int_0^T\xi_u\;{\rm
d}S_u\right|\cG_0\right)}{\IE_{\IQ_\IG}H}=\alpha.\nonumber \een
Applying second part of Lemma \ref{NPL} for the success set
$A=\left\{\alpha\IE_{\IQ_\IG}H+\int_0^T\xi_u\;{\rm d}S_u\geq
H\right\}$ and using required inequality (\ref{P2inequality}) and
definition of $k$ given in (\ref{eqk2}) we derive: \ben
\alpha&\geq&\IQ^*\left(\left.\alpha\IE_{\IQ_\IG}H+\int_0^T\xi_u\;{\rm d}S_u\geq H\right|\cG_0\right)\nonumber\\
&\geq&\IQ^*\left(\left.\left.\frac{{\rm d}\IQ^*}{{\rm d}\IP}\right|_{\cG_T}\leq k\right|\cG_0\right).\nonumber
\label{rhs}
\een
We prove now that for this particular minimal choice of $\alpha$ being the rhs of (\ref{rhs})
the strategy $\left(\IE_{\IQ_\IG}\left[H \left.\mathbf{1}_{\left\{\left.\frac{{\rm d}\IQ^*}{{\rm d}\IP}\right|_{\cG_T}\leq k\right\}}\right|\cG_0\right],\tilde{\xi}\right)$ satisfies the inequality (\ref{P2inequality}) of the Problem \ref{P2}:
\ben
\lefteqn{\IP\left(\left.\alpha\IE_{\IQ_\IG}H+\int_0^T\tilde{\xi}_u\;{\rm d}S_u\geq H\right|\cG_0\right)}\nonumber\\
&&=\IP\left(\left.\IE_{\IQ_\IG}\left.\left[H \mathbf{1}_{\left\{\left.\frac{{\rm d}\IQ^*}{{\rm d}\IP}\right|_{\cG_T}\leq k\right\}}\right|\cG_0\right]
+\int_0^T\tilde{\xi}_u\;{\rm d}S_u\geq H\right|\cG_0
\right)\nonumber\\
&&=
\IP\left(\left. H \mathbf{1}_{\left\{\left.\frac{{\rm d}\IQ^*}{{\rm d}\IP}\right|_{\cG_T}\leq k\right\}}
\geq H\right|\cG_0
\right)
\nonumber\\
&&\geq\IP\left(\left.\left.\frac{{\rm d}\IQ^*}{{\rm d}\IP}\right|_{\cG_T}\leq k\right|\cG_0\right)=1-\epsilon.\nonumber
\een
This completes the proof.
\halmos

\section*{Acknowledgements}
This work is partially supported by the Ministry of Science and Higher Education of Poland under the grant N N2014079 33 (2007-2009).

\bibliographystyle{abbrv}


\end{document}